\documentclass[11pt]{amsart}

\usepackage{amsmath,amstext,amscd}
\usepackage{amssymb}
\usepackage{epsfig}
\usepackage{amsfonts,latexsym}
\usepackage{hyperref}

\newtheorem{defin}{Definition}[section]
\newtheorem{thm}[defin]{Theorem}
\newtheorem{cor}[defin]{Corollary}
\newtheorem{rem}[defin]{Remark}
\newtheorem{lem}[defin]{Lemma}
\newtheorem{prop}[defin]{Proposition}
\newtheorem{quest}[defin]{Algorithmn}
\newtheorem{prob}[defin]{Notation}

\newcommand{\theorem}[1]{\begin{thm} \sl{#1} \end{thm}}
\newcommand{\theoremname}[2]{\begin{thm}[#1] \sl{#2} \end{thm}}

\newcommand{\lemma}[1]{\begin{lem} \sl{#1} \end{lem}}

\newcommand{\remark}[1]{\begin{rem} \emph{#1} \end{rem}}

\newcommand{\corollary}[1]{\begin{cor} \sl{#1} \end{cor}}
\newcommand{\corollaryname}[2]{\begin{cor}[#1] \sl{#2} \end{cor}}
\newcommand{\lemmaname}[2]{\begin{lem}[#1] \sl{#2} \end{lem}}
\newcommand{\problem}[1]{\begin{prob} \emph{#1} \end{prob}}

\begin{document}

\title{Cryptanalysis of the Shpilrain-Ushakov protocol for Thompson's group}

\author{Francesco Matucci}

\dedicatory{\emph{Department of Mathematics, Cornell University, Ithaca, NY 14853, USA} \\
\tt{matucci@math.cornell.edu}}

\begin{abstract}
This paper shows that an eavesdropper can always recover efficiently the private key
of one of the two parts of the public key cryptography protocol introduced by Shpilrain
and Ushakov in \cite{su}. Thus an eavesdropper can always recover the shared secret key,
making the protocol insecure. 
\end{abstract}

\maketitle

\markboth{The Shpilrain-Ushakov Protocol for Thompson's Group $F$ is always breakable}{Francesco Matucci}

\section{Introduction}
Recent advances in public key cryptography have underlined the need to find alternatives
to the RSA cryptosystem. It has been proposed to use algorithmic problems
in non-commutative group theory as possible ways to build new protocols. 
The \emph{conjugacy search problem} was introduced in several papers
as a generalization of the \emph{discrete logarithm problem} in the research of a new safe encryption scheme.
The former problem asks whether or not, given a group $G$ and two elements $a,b \in G$ that are conjugate, 
we can find at least one $x \in G$ with $a^x:=x^{-1}ax=b$. It is thus important to look for a platform
group $G$ where this problem is computationally hard. Seminal works by Anshel-Anshel-Godlfeld \cite{aag} and Ko-Lee et al.\ \cite{kolee}
have proposed the braid group $B_n$ on $n$ strands as a possible platform group.

\medskip

It has been observed that Thompson's group $F$ and the braid groups $B_n$ have some similarities.
Belk proved in his thesis \cite{belkthesis} that $F$ and the braid groups have a similar classifying
space. Loosely speaking, the elements of $F$ appear as braids, but with merges and splits
instead of twists (this representation of $F$ uses \emph{strand diagrams} which are introduced in \cite{belkthesis}).
Dehornoy defined in \cite{deho} a group of \emph{parenthesized braids} which contains both $F$ and $B_n$ in
a very natural way. However, for cryptographic purposes, $F$ has still not proved to be a
good platform. Kassabov and Matucci have proved in \cite{kama}
that the simultaneous conjugacy problem is efficiently solvable, making it insecure to apply
protocols based on the conjugacy problem.

\medskip

Shpilrain and Ushakov in \cite{su} have proposed using a particular version of the
\emph{decomposition problem} as a protocol and the group $F$ as a platform. The new problem is: given a
group $G$, a subset $X \subseteq G$ and two elements $w_1,w_2 \in G$ with the information that there exist
$a,b \in X$ such that $a w_1 b =w_2$, find at least one such pair $a,b$.
In this paper we show how to recover efficiently the shared secret key of this protocol.

\medskip

The paper is organized as follows. In Section \ref{sec:protocol} and Section \ref{sec:platform}
we recall the protocol and give a description of Thompson's group $F$. In Section 4 we
recall the choice of parameters proposed in \cite{su}.
In section \ref{sec:recovering-keys} we give an efficient attack that always recovers
the secret key. In Sections \ref{sec:transitivity} and
\ref{sec:recovering-keys-alternatives} we show another type of attack.
In Section \ref{sec:restriction-centralizers} we make some comments on possible generalizations of this protocol.

\subsection*{History and related works.} The first attack on this protocol
was announced by Ruinskiy, Shamir and Tsaban in November 2005 at the Bochum Workshop 
\emph{Algebraic Methods in Cryptography}, showing that the paramaters given in \cite{su}
should be increased to have higher security of the system. Their attack was improved in other
announcements and was finalized in \cite{rst2} at the same time that this paper was written.
Their attack describes a more general procedure which uses length functions.
We remark that the same authors have been developing new techniques involving
``subgroup distance functions'' and that they applied them on the same protocol for $F$ as a test case \cite{rst3}.
The approach of Ruinskiy, Shamir and Tsaban in their mentioned papers is heuristic,
and its success rates are good but not $100\%$. Our approach is deterministic, and provably succeeds
in all possible cases.

\subsection*{Acknowledgements.} The author would like to thank Martin Kassabov, Boaz Tsaban and Vladimir Shpilrain
for helpful discussions. The author would also like to thank Ken Brown and the referees for many helpful comments.

\section{The Protocol \label{sec:protocol}}

The protocol proposed in \cite{su} is based on the \emph{decomposition problem}:
given a group $G$, a subset $X \subseteq G$ and $w_1,w_2 \in G$, find $a,b \in X$ with $aw_1b=w_2$, given that such $a,b$ exist. Here is
the protocol in detail:

\subsection*{Public Data.} A group $G$, an element $w \in G$ and two subgroups $A,B$ of $G$
such that $ab=ba$ for all $a \in A$, $b \in B$.

\subsection*{Private Keys.} Alice chooses $a_1 \in A$, $b_1 \in B$ and sends the element
$u_1=a_1 w b_1$ to Bob. Bob chooses $b_2 \in B$, $a_2 \in A$ and sends the element $u_2=b_2wa_2$ to Alice.
Alice then computes the element $K_A=a_1 u_2 b_1=a_1 b_2 w a_2 b_1$ and Bob computes the element
$K_B=b_2 u_1 a_2=b_2 a_1 w b_1 a_2$. Since $A$ and $B$ commute elementwise, $K=K_A=K_B$ becomes
Alice and Bob's shared secret key.

\subsection*{Eavesdropper's Data.} Eve has all the public data and the two elements $u_1$ and
$u_2$, observed during Alice and Bob's exchange.

\section{The Group $F$ and the Subgroups $A_s,B_s$ \label{sec:platform}}

Thompson's group $F$ was introduced by R. Thompson while working on problems in logic.
The standard introduction to $F$ is \cite{cfp}. One of Thompson's original definitions of $F$ is the following:
for $I=[0,1]$ we define $PL_2(I)$ to be the group
of piecewise linear homeomorphisms of the interval $I$
with finitely many breakpoints such that:
\begin{itemize}
\item
all slopes are integral powers of $2$, and
\item
all breakpoints are in $\mathbb{Z}[\frac{1}{2}]$, the ring of dyadic rational numbers;
\end{itemize}
the product of two elements is given by the composition of functions. We thus define $F$ to be the group
$PL_2(I)$. $F$ can also be described using the following presentation:
\[
F = \langle x_0, x_1, x_2, \ldots \mid x_n x_k = x_k x_{n+1}, \forall \, k<n\rangle.
\]
This presentation has the advantage that the elements of $F$ can be uniquely written in the following \emph{normal form}
\[
x_{i_1}\ldots x_{i_u} x_{j_v}^{-1} \ldots x_{j_1}^{-1}
\]
such that $i_1 \le \ldots \le i_u$, $j_1 \le \ldots \le j_v$ and if both $x_i$ and
$x_i^{-1}$ occur, then either $x_{i+1}$ or $x_{i+1}^{-1}$ occurs, too. Since $x_k=x_0^{1-k}x_1x_0^{k-1}$
for $k \ge 2$, the group $F$ is generated by the elements $x_0$ and $x_1$.
The generators $x_k$ of the infinite presentation can be represented as piecewise-linear homeomorphisms
by shrinking the function $x_0$ shown in figure \ref{fig:generators-F} onto the interval
$[1-\frac{1}{2^k},1]$ and extending it as the identity on $[0,1-\frac{1}{2^k}]$:

\begin{figure}[0.5\textwidth]
 \begin{center}
  \includegraphics[height=6cm]{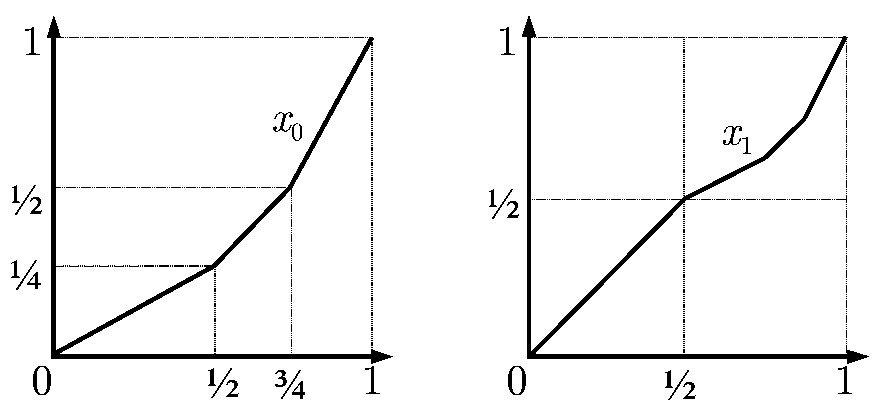}
 \end{center}
 \caption{Two of the elements of the generating set of $F$.}
 \label{fig:generators-F}
\end{figure}

\begin{figure}[0.5\textwidth]
 \begin{center}
  \includegraphics[height=6cm]{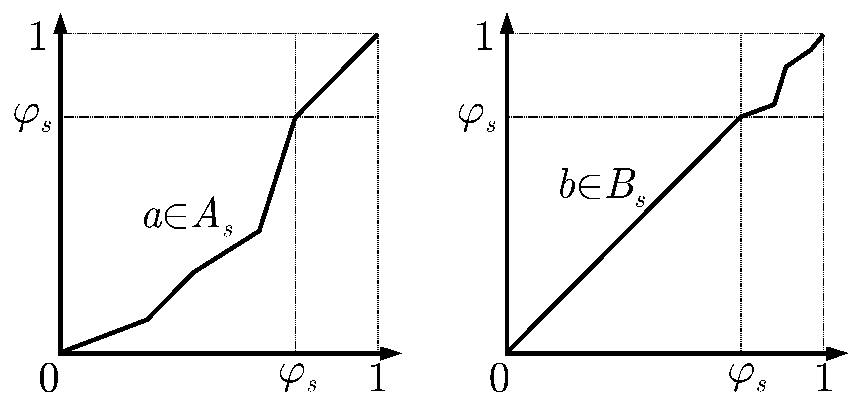}
 \end{center}
 \caption{An example of an element of $A_s$ and one of $B_s$.}
 \label{fig:commuting-elements}
\end{figure}

\noindent We now introduce a notation which will be useful for the definition of the subgroups $A$ and $B$. 
For every positive integer $k$ we call
\[
\varphi_k:=1-\frac{1}{2^{k+1}}.
\]
From the definition of $x_k$, we get
\[
x_k^{-1}\left(\left[\varphi_k,1\right]\right)=\left[\varphi_{k+1},1\right] \subseteq
\left[\frac{3}{4},1\right]
\]
implying that, for $t \in [\varphi_k,1]$, we have
\[
\frac{d}{dt}x_0x_k^{-1}(t)=x_0'(x_k^{-1}(t))(x_k^{-1})'(t)=2 \cdot \frac{1}{2}=1
\]
which means $x_0 x_k^{-1}$ is the identity in the interval $[\varphi_k,1]$.
For any $s \in \mathbb{N}$, Shpilrain and Ushakov define in \cite{su} the following sets
\[
S_{A_s}=\{x_0 x_1^{-1}, \ldots, x_0 x_s^{-1}\}
\]
\noindent and
\[
S_{B_s}=\{x_{s+1},\ldots,x_{2s}\}
\]
and then define the subgroups $A_s:=\langle S_{A_s} \rangle$ and $B_s:=\langle S_{B_s} \rangle$.
The previous argument immediately yields that that all elements of $A_s$ commute with
all elements of $B_s$ (see figure \ref{fig:commuting-elements}), i.e.

\lemmaname{Shpilrain-Ushakov \cite{su}}{For every fixed $s \in \mathbb{N}$, $ab=ba$ for every
elements $a \in A_s$ and $b \in B_s$.}

\problem{For every dyadic number $d \in [0,1]$ we denote by $PL_2([0,d])$ the set of functions in $PL_2(I)$
which are the identity on $[d,1]$. Moreover, if we are given a
piecewise linear map defined only on $[0,d]$ we will assume it is extended to $[0,1]$ by defining it
as the identity on $[d,1]$. Similar remarks apply to $PL_2([d,1])$.}

\noindent Parts (i) and (iii) of the following Lemma are in \cite{su}, while part (ii) is a simple observation.

\lemma{(i) $A_s$ is the set of elements whose normal form is of the type
\[
x_{i_1}\ldots x_{i_m} x_{j_m}^{-1} \ldots x_{j_1}^{-1}
\]
where $i_k - k < s$ and $j_k - k < s$, for all $k=1,\ldots,m$.

\medskip
\noindent (ii) $B_s=PL_2([\varphi_s,1])$.

\medskip
\noindent (iii) Let $a \in A_s$ and $b \in B_s$ be such that their normal forms are
\begin{eqnarray*}
a = x_{i_1}\ldots x_{i_m} x_{j_m}^{-1} \ldots x_{j_1}^{-1} \\
b = x_{c_1}\ldots x_{c_u} x_{d_v}^{-1} \ldots x_{d_1}^{-1}.
\end{eqnarray*}
Then the normal form of $ab$ is
\[
ab = x_{i_1}\ldots x_{i_m} x_{c_1+m}\ldots x_{c_u+m} x_{d_v+m}^{-1} \ldots x_{d_1+m}^{-1}
x_{j_m}^{-1} \ldots x_{j_1}^{-1}.
\]\label{thm:definition-A_s}}

\theoremname{Shpilrain-Ushakov \cite{su}}{In Thompson's group $F$, the normal form of a given word $w$
can be computed in time $O(|w| \log |w|)$,
where $|w|$ is the length of the normal form in the generators $x_0,x_1, x_2, \ldots$ \label{thm:normal-forms}}

\section{Suggested Parameters for the Encryption \label{sec:parameters}}

We now illustrate briefly the choice of parameters proposed in \cite{su}.
Alice and Bob select an integer $s \in [3,8]$ and an even integer
$M \in [256,320]$ uniformly and randomly.
Morever, they also choose a random element $w \in \langle x_0,x_1, \ldots, x_{s+2} \rangle$ with $|w|=M$,
where $|w|$ is as in Theorem \ref{thm:normal-forms}.
The numbers $s,M$ and the element $w$ are now part of the the public data.

To proceed with the protocol described in Section \ref{sec:protocol}, 
Alice chooses random elements $a_1 \in A_s,b_1 \in B_s$, with $|a_1|=|b_1|=M$, while Bob chooses
random elements $a_2 \in A_s,b_2 \in B_s$, with $|a_2|=|b_2|=M$.
Now they both compute the shared secret key:
\[
K=a_1 b_2 w a_2 b_1.
\]
Shpilrain and Ushakov remark that this choice of parameters gives a key space which 
increases exponentially in $M$, i.e., $|A_s(M)| \ge \sqrt{2}^M$, thereby making it difficult for Eve to perform a brute force attack.

\section{Recovering the Shared Secret Key \label{sec:recovering-keys}}

We begin this section by providing the theoretical background
for the attack. We will use the piecewise-linear point of view to understand why the attack works
and then rephrase it combinatorially.
We will now describe how Eve, by knowing the elements $w,u_1,u_2$, can always
recover one of the two legitimate parties' private keys. She chooses whose key to crack, depending on whether the graph of $w$
is above or below the point $(\varphi_s,\varphi_s)$. 

\subsection{Recovering Bob's Private Keys: $w(\varphi_s) \le \varphi_s$}

Since $w(t)\le \varphi_s$ for all $t \in [0,\varphi_s]$, we observe the following identity
\[
u_2(t)=b_2 w a_2(t)=w a_2(t), \, \, \forall t \in [0,\varphi_s].
\]
Therefore, Eve may apply $w^{-1}$ to the left of both sides of the previous equation to obtain
\[
w^{-1} u_2(t)=a_2(t), \, \, \forall t \in [0,\varphi_s]
\]
and so $w^{-1}u_2 \in A_s B_s$ and
\[
a_2(t)=\begin{cases}
w^{-1} u_2(t) & t \in [0,\varphi_s] \\
t & t \in [\varphi_s,1].
\end{cases}
\]
Now Eve has the elements $a_2$, $w$ and $u_2=b_2wa_2$ and she computes
\[
b_2=u_2 a_2^{-1}w^{-1}
\]
thereby detecting Bob's private keys and the shared secret key $K$.

\subsection{Recovering Alice's Private Key: $w(\varphi_s)>\varphi_s$}

Since $w^{-1}(t) < \varphi_s$ for all $t \in [0,\varphi_s]$, we have
\[
u_1^{-1}(t)=b_1^{-1} w^{-1} a_1^{-1}(t)= w^{-1} a_1^{-1}(t)  ,\, \, \forall t \in [0,\varphi_s].
\]
By applying the same technique as in the previous subsection Eve recovers $a_1^{-1}$ and obtains that $u_1w^{-1} \in A_s B_s$. Thus,
she is able to detect $a_1,b_1$ and the shared secret key $K$. Alternatively, Eve observes
\[
w^{-1}u_1(t)=w^{-1}a_1 w b_1(t)= b_1(t)  ,\, \, \forall t \in [\varphi_s,1]
\]
and so
\[
b_1(t)=\begin{cases}
t & t \in [0,\varphi_s] \\
w^{-1} u_1(t) & t \in [\varphi_s,1].
\end{cases}
\]

\subsection{Outline of the attack} We expand on the previous discussion to describe a combinatorial attack.
Assume that Eve has the elements $w,u_1,u_2$.

\begin{enumerate}
\item{Eve writes the normal forms of $z_1:=u_1 w^{-1}$ and $z_2:=w^{-1}u_2$.}
\item{By the previous discussion, either $z_1 \in A_s B_s$ or $z_2 \in A_s B_s$ (or both).
She can detect which one using Lemma \ref{thm:definition-A_s}(i) and selects this $z_i$.}
\item{She  computes the $A_s$-part $a_{z_i}$ of $z_i$.}
\item{If $i=1$, she computes $b_{z_1}:= w^{-1} a_{z_1}^{-1} u_1$. If $i=2$,
she computes $b_{z_2} := u_2 a_{z_2}^{-1}w^{-1}$.}
\item{Eve computes $K$ from $u_1,u_2,a_{z_{i}},b_{z_i}$.}
\end{enumerate}

The only point of this procedure which needs further explanation is (2). When we have
the normal forms of $z_1,z_2$, we know that one of them is in $A_s B_s$.
We write the normal form $z_i = x_{i_1}\ldots x_{i_e} x_{j_f}^{-1} \ldots x_{j_1}^{-1}$
and we look at the notation of Lemma \ref{thm:definition-A_s}(i): we need to find
the the smallest index $r$ in $z_i$ such that either $i_{r+1}$ or $j_{r+1}$
does not satisfy the index condition in Lemma \ref{thm:definition-A_s}(i). To verify if $z_i \in A_s B_s$,
we need to check whether it has the form described in Lemma \ref{thm:definition-A_s}(iii): we remove the first
$r$ letters and the last $r$ letters of $z_i$ from the word and we lower all the indices of the remaining letters by $r$;
if what remains is a word whose indices are in $\{s+2, \ldots, 2s\}$, then we have an element of $B_s$, otherwise $z_i \not \in A_s B_s$.
If $z_i \in A_s B_s$, then $a_{z_i}$ will be 
the product of the first $r$ elements of $z_i$ and the last $r$ ones.

\subsection{Complexity of the attack.} By Theorem \ref{thm:normal-forms}
we know that computing normal forms can be done in time $O(M \log M)$, where $M$
is the size of the inputs suggested in Section \ref{sec:parameters}. Part (2) of the attack can be executed in time $O(M)$,
by just reading the indices of the normal forms and finding when the relation of Lemma
\ref{thm:definition-A_s}(i) breaks down. Finally, the last steps are just
multiplications and then simplifications so they can again be performed in time
$O(M \log M)$. Therefore, Eve can recover the shared secret key in time $O(M \log M)$.

\remark{The previous discussion shows that there is no need to pass from words to piecewise-linear functions and back.
The attack can be performed entirely by using the combinatorial point of view which is used for encryption. The piecewise-linear
point of view is necessary only to prove that the combinatorial attack works. We also remark that the complexity of
the attack is independent of the parameter $s$.}

\section{Transitivity of $A_s$ and $B_s$ \label{sec:transitivity}}

The previous section showed how to recover the shared secret key of one of
the two involved parties, based on whether the graph of $w$ lies above or below
the point $(\varphi_s,\varphi_s)$. However, it is possible to find the shared
secret key even in the cases not studied in the previous section. More precisely,
it is possible to attack Alice's word in the case $w(\varphi_s) \le \varphi_s$
and Bob's word in the case $w(\varphi_s)>\varphi_s$. We need a better description of
the subroups $A_s$. If $s=1$, we observe that $A_1=\langle x_0 x_1^{-1}\rangle$ is a cyclic group.
For larger values of $s$, $A_s$ becomes the full group of piecewise linear homeomorphism
on $[0,\varphi_s]$.

\lemma{$A_2=PL_2\left(\left[0,\frac{7}{8}\right]\right)$. \label{thm:first-isomorphism}}

\begin{figure}[0.5\textwidth]
 \begin{center}
  \includegraphics[height=6cm]{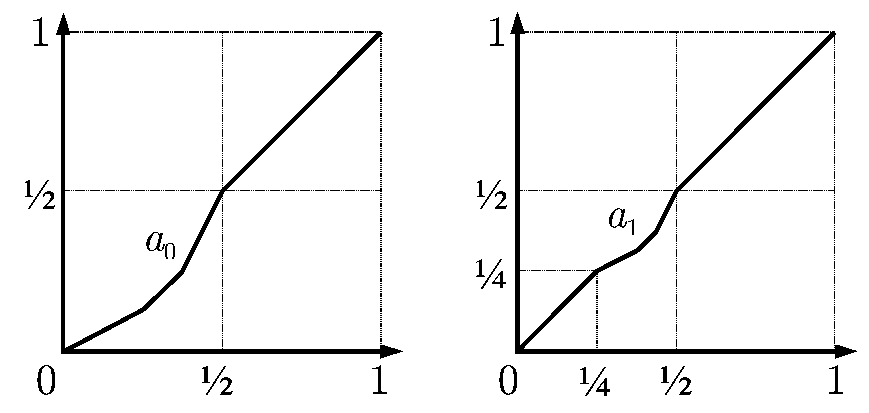}
 \end{center}
 \caption{The two standard generators for $PL_2([0,\frac{1}{2}])$.}
 \label{fig:generators-A_s}
\end{figure}

\noindent \emph{Proof.} Let $a,b$ be the two generators of $PL_2([0,\frac{1}{2}])$ shown
in figure \ref{fig:generators-A_s}.

\noindent One sees that $a=x_0^2x_1^{-1}x_0^{-1}$
and that $b=x_0 x_1^2 x_2^{-1} x_1^{-1} x_0^{-1}$ and so a conjugation of $PL_2([0,\frac{1}{2}])$
by $x_0^2$ yields $PL_2([0,\frac{7}{8}])=\langle x_0^2 a x_0^{-2}, x_0^2 b x_0^{-2}\rangle$. By Lemma
\ref{thm:definition-A_s} we have
\begin{eqnarray*}
x_0^2 a x_0^{-2} = x_0^4 x_1^{-1}x_0^{-3} \in A_2 \\
x_0^2 b x_0^{-2} = x_0^3 x_1^2 x_2^{-1} x_1^{-1} x_0^{-3} \in A_2
\end{eqnarray*}
so that $PL_2([0,\frac{7}{8}]) \subseteq A_2$. The other inclusion is obvious. $\square$

\theorem{$A_s=PL_2([0,\varphi_s])$, for every $s \ge 2$. \label{thm:subgroups-are-F}}

\noindent \emph{Proof.} A straightforward computation shows that
\[
x_0^{-1} PL_2([0,\varphi_s])x_0=PL_2([0,\varphi_{s+1}]), \forall s\ge0.
\]
Therefore $A_2=PL_2([0,\varphi_2])$ and the definition of $A_s$ imply

\[
PL_2([0,\varphi_s])= x_0^{s-2} A_2 x_0^{2-s} \subseteq A_s \subseteq PL_2([0,\varphi_s])
\]

\medskip

\noindent therefore implying that $A_s=PL_2([0,\varphi_s])$. $\square$

\corollary{$A_s \cong B_s \cong F$, for every $s \ge 2$.}

\noindent The previous Theorem and Lemma 2.5 in \cite{kama} yield the following corollaries:

\corollaryname{Transitivity of $A_s$}{For any $t_1,t_2 \in \mathbb{Z}\left[\frac{1}{2}\right] \cap [0,\varphi_s]$ we can construct an
$a \in A_s$ with $a(t_1)=t_2$. \label{thm:transitivity-A_s}}

\corollaryname{Extendability of $A_s$}{Let $t_0 \in \mathbb{Z}\left[\frac{1}{2}\right] \cap [0,\varphi_s]$ and $\bar{a}(t)=a|_{[0,t_0]}$
for an element $a \in A_s$. Assume we know $\bar{a}$, but that we do not know $a$. Then we can construct
an $a_\sigma \in A_s$ such that $a_\sigma(t)=\bar{a}(t)$ for all $t \in [0,\varphi_s]$.
\label{thm:extension-A_s}}

\remark{The analogues of the last two corollaries are true for the interval
$[\varphi_s,1]$ and $B_s$ too.}

\section{Using Transitivity to Attack the Shared Secret Key \label{sec:recovering-keys-alternatives}}

With the new description of $A_s$ and $B_s$ given in
section \ref{sec:transitivity}, it is now possible to attack the secret
keys in the cases left open from section \ref{sec:recovering-keys}.

\subsection{Attacking Alice's word for the case $w(\varphi_s) \le \varphi_s$}

We have
\begin{eqnarray*}
u_1(t)=a_1w(t), \forall t \in [0,\varphi_s],
\end{eqnarray*}
thus
\begin{eqnarray*}
a_1(t)=u_1w^{-1}(t), \forall t \in [0,w(\varphi_s)]
\end{eqnarray*}
and so $a_1$ is uniquely determined in $[0,w(\varphi_s)]$.
We apply corollary \ref{thm:extension-A_s} to find an element $a_\sigma \in A_s$ such
that $a_\sigma=a_1$ on the interval $[0,w(\varphi_s)]$. If we define
\[
b_\sigma:=w^{-1} a_\sigma^{-1} u_1
\]
then we have that
\begin{eqnarray*}
b_\sigma(t)=w^{-1}a_\sigma^{-1}a_1 w(t)= w^{-1} w(t)=t, \forall t \in [0,\varphi_s]
\end{eqnarray*}

\noindent Therefore $b_\sigma \in B_s$ and $a_\sigma w b_\sigma = u_1$
and so Eve can recover the shared secret key $K$ by using the pair $(a_\sigma,b_\sigma)$.

\remark{We observe
that any extension of $a_1|_{[0,w(\varphi_s)]}$ to an element $a_\sigma$ of $PL_2([0,\varphi_s])$ will yield a suitable element
to attack Alice's key. Moreover, any element $a_1' \in A_s$
such that $a_1'wb_1'=u_1$, for some suitable $b_1' \in B_s$,
will be an extension of $a_1|_{[0,w(\varphi_s)]}$.}

\subsection{Attacking Bob's word for the case $w(\varphi_s) > \varphi_s$}

Eve considers $u_2^{-1}=a_2^{-1} w^{-1} b_2^{-1}$ and recovers a pair $(a_\sigma^{-1},b_\sigma^{-1})$
to get the shared secret key in the same fashion of the previous subsection.

\remark{Both the techniques of this section have been carried out using the transitivity of $A_s$ (Corollary \ref{thm:transitivity-A_s}).
They can also be solved by using the analogue of Corollary
\ref{thm:extension-A_s} for $B_s$ to get another pair $(a_\sigma,b_\sigma)$ which can be used to
retrieve the secret key.}

\section{Comments and Alternatives to the Protocol \label{sec:restriction-centralizers}}

This section analyzes possible alternatives and weaknesses of our methods. We observe that,
if instead of $PL_2(I)$ we had used a larger group of piecewise linear homeomorphisms of the unit interval, the same
technique would have worked, as long as the commuting subgroups $A$ and $B$ had disjoint supports.
More generally, we can copy this idea if the given group $G$ acts on some space and we
have $A, B$ with disjoint support. We will now see some examples
of how this is possible.

\subsection{Choice of the subgroups $A$ and $B$}

We recall the following result:

\theoremname{Kassabov-Matucci \cite{kama}}{Let $A=\langle a_1, \ldots, a_m \rangle \le F$ be a finitely generated subgroup. Then

\medskip
\noindent (i) There exists a dyadic partition of $[0,1]=I_1 \cup \ldots \cup I_n$ such that the centralizer
$C_F(A):=\{f \in F \, | \, af=fa, \forall a \in A\}$ 
is a product of subgroups $C_1, \ldots, C_n$, where $C_r \le \{f \in F \, | \, f(t)=t, \forall t \not \in I_r \}$. 
Moreover, we have
\begin{itemize}
\item $C_r=PL_2(I_r)$ if and only if of $a_i|_{I_r}=id$, for all $i=1, \ldots, r$.
\item $C_r \cong \mathbb{Z}$ if and only if $a_1|_{I_r}, \ldots, a_m|_{I_r}$ have a common root on $I_r$.
\item $C_r = 1$ if and only if there are $i \ne j$ such that $a_i|_{I_r}, a_j|_{I_r}$ have no common root on $I_r$.
\end{itemize}

\medskip
\noindent (ii) There exist two elements $g_1,g_2 \in F$ such that $C_F(A)=C_F(g_1) \cap C_F(g_2)$.}

Going back to the protocol introduced in Section \ref{sec:protocol} we observe that, after
we choose a finitely generated subgroup $A=\langle f_1,\ldots,f_m \rangle$,
we are very restricted in our choice of the subgroup $B$. Since $B \le C_F(A)$, we must make sure that the elements of
$B$, when restricted to $I_r$, are powers of common roots of the $a_i$'s, if at least one $a_i$ is non-trivial on $I_r$.
This gives a tight restriction on the subgroup $B$ whose support is essentially disjoint from that of $A$, except in the
intervals where they all are powers of a common root. An attack similar to that of Section \ref{sec:recovering-keys}
can thus be applied on each interval $I_r$: if their supports are disjoint on $I_r$, we can act as before, otherwise elements
of $A$ and $B$ are powers of a common root on $I_r$.

With more general commuting subgroups, the attack of Section \ref{sec:recovering-keys} does not immediately give either of the two keys.
However, the discussion above suggests that the choice of $A$ and $B$ must be done much more carefully in order to avoid similar
attacks.

\medskip

\subsection{Alternative Protocol and Attacks}

Ko-Lee et al.\ \cite{kolee} introduced a slightly different protocol based on the decomposition problem
(They worked with braid groups, but we will apply their protocol to Thompson's group).
In their protocol, Alice picks $a_1,a_2 \in A$ and sends $u_1 = a_1 w a_2$ to Bob, while Bob
chooses $b_1,b_2 \in B$ and sends $u_2=b_1 w b_2$ to Alice. We can still attempt to solve this new protocol, by
again dividing the problem into various cases. We assume that we use the same subgroups $A_s$ and $B_s$
and we work in the case $w(\varphi_s) \le \varphi_s$ to
show how to attack the private keys of Bob.
We apply the analogue for $B_s$ of Corollary \ref{thm:transitivity-A_s} and
find a $b_0$ such that $b_0^{-1}(w^{-1}(\varphi_s))=u_2^{-1}(\varphi_s)= b_2^{-1}w^{-1}(\varphi_s)$.
We define
\begin{eqnarray*}
b_1'=b_1 \\
b_2'= b_2 b_0^{-1} \\
u_2'=b_1'wb_2'
\end{eqnarray*}
so that $b_2'(w^{-1}(\varphi_s))=w^{-1}(\varphi_s)>\varphi_s$. Thus we have
\[
u_2'(t)=b_1'(t)wb_2'(t)=wb_2'(t), \forall t \in [0,w^{-1}(\varphi_s)]
\]
hence
\[
b_2'(t)=w^{-1}u_2'(t), \forall t \in [0,w^{-1}(\varphi_s)].
\]

\noindent Thus $b_2'$ is uniquely determined in $[0,w^{-1}(\varphi_s)]$.
We apply corollary \ref{thm:extension-A_s} for $B_s$ to find
a $b_{\sigma_2} \in B_s$ such that $b_{\sigma_2}=b_2'$ on $[0,w^{-1}(\varphi_s)]$
and we define
\[
b_{\sigma_1}:=u_2'b_{\sigma_2}^{-1}w^{-1}.
\]
Thus
\[
b_{\sigma_1}(t)=b_1' w b_2' b_{\sigma_2}^{-1}w^{-1}(t)=b_1'(t)=t, \forall t \in [0,\varphi_s]
\]
therefore $b_{\sigma_1} \in B_s$. Therefore the pair $(b_{\sigma_1},b_{\sigma_2})$ satisfies
$u_2'=b_{\sigma_1}wb_{\sigma_2}$ and so Eve can recover the shared secret key $K$. A similar argument
can be used to attack the element $a_1wa_2$, with the transitivity results for $A_s$.

\subsection{A comment on the Alternative Protocol}

The weakness in the protocol discussed in the previous
subsection arises from the fact that the chosen subgroups $A_s$ and $B_s$ are transitive on the intervals
on which they act nontrivially. This suggests that a possible way to
avoid such attacks is for $A$ and $B$ to be chosen to be not transitive on their support.

\remark{We observe that the attacks of section \ref{sec:recovering-keys-alternatives} and section
\ref{sec:restriction-centralizers} can be carried out in a fashion similar to that of Section \ref{sec:recovering-keys}, still producing
a solution in polynomial time.}

\bibliographystyle{plain}

\begin{thebibliography}{1}

\bibitem{aag}
I.~Anshel, M.~Anshel, and D.~Goldfeld.
\newblock An algebraic method for public-key cryptography.
\newblock {\em Math. Res. Lett.}, 6(3-4):287--291, 1999.

\bibitem{belkthesis}
J.M. Belk.
\newblock {\em Thompson's Group {$F$}}.
\newblock PhD thesis, Cornell University, 2004.
\newblock \url{http://arxiv.org/abs/0708.3609}

\bibitem{cfp}
J.W. Cannon, W.J. Floyd, and W.R. Parry.
\newblock Introductory notes on {R}ichard {T}hompson's groups.
\newblock {\em Enseign. Math. (2)}, 42(3-4):215--256, 1996.

\bibitem{deho}
P.~Dehornoy.
\newblock The group of parenthesized braids.
\newblock {\em Adv. Math.}, 205(2):354--409, 2006.
\newblock \url{http://arxiv.org/abs/math/0407097}

\bibitem{kama}
M.~Kassabov and F.~Matucci.
\newblock The simultaneous conjugacy problem in {T}hompson's group {$F$}.
\newblock {\em Groups, Geometry and Dynamics}, to appear
\newblock \url{http://arxiv.org/abs/math/0607167}

\bibitem{kolee}
K.H. Ko, S.J. Lee, J.H. Cheon, J.W. Han, J.~Kang, and C.~Park.
\newblock New public-key cryptosystem using braid groups.
\newblock In {\em Advances in cryptology---CRYPTO 2000 (Santa Barbara, CA)},
  volume 1880 of {\em Lecture Notes in Comput. Sci.}, pages 166--183. Springer,
  Berlin, 2000.

\bibitem{rst2}
D.~Ruinskiy, A.~Shamir, and B.~Tsaban.
\newblock Length-based cryptanalysis: The case of {T}hompson's group.
\newblock {\em Journal of Mathematical Cryptology} \textbf{1}, 2007,
no. 4, 359--372, \url{http://arxiv.org/abs/cs/0607079}

\bibitem{rst3}
D.~Ruinskiy, A.~Shamir, and B.~Tsaban.
\newblock Cryptanalysis of group-based key agreement protocols using subgroup
  distance functions.
\newblock In {\em Proceedings of the 10th International Conference on Practice
  and Theory in Public-Key Cryptography PKC07}, volume 4450 of {\em Lecture
  Notes in Comput. Sci.}, pages 61--75. 2007. \url{http://arxiv.org/abs/0705.2862}

\bibitem{su}
V~Shpilrain and A.~Ushakov.
\newblock Thompson's group and public key cryptography.
\newblock In {\em ACNS 2005}, volume 3531 of {\em Lecture Notes in Comput.
  Sci.}, pages 151--163. 2005.
\newblock \url{http://arxiv.org/abs/math/0505487}

\end{thebibliography}

\end{document}